\documentclass[11pt]{article}
\usepackage[affil-it]{authblk}
\usepackage{natbib}
\bibpunct[, ]{(}{)}{,}{a}{}{,}%
\usepackage[font=footnotesize]{caption}
\usepackage{color}
\usepackage{amssymb}
\usepackage{amsmath}
\usepackage{amsthm}
\usepackage[margin=1in]{geometry}
\usepackage{setspace}
\usepackage{mathtools}
\usepackage{titlesec}
\usepackage{algorithm}
\usepackage{algpseudocode}
\usepackage{enumitem}
\usepackage{array}
\usepackage{tabularx}
\usepackage{multirow}

\usepackage[titletoc,toc,title]{appendix}
\usepackage[pagewise]{lineno}

\titleformat{\section}                                     
{\sffamily\bfseries}
{\thesection.}{1em}{}
\titleformat{\subsection}                                   
{\sffamily\slshape\bfseries}
{\thesubsection.}{1em}{}

\title{\singlespacing{Robust Binary Linear Programming Under Implementation Uncertainty}}
\author[1]{Jose E. Ramirez-Calderon}
\affil[1]{\small{Department of Industrial and Systems Engineering, Texas A\&M University,3131 TAMU, College Station, TX 77843-3131, USA, ramirez.jose@tamu.edu}}
\author[2]{V. Jorge Leon}
\affil[2]{\small{Department of Engineering Technology and Industrial Distribution and Department of Industrial and Systems Engineering, Texas A\&M University,3367 TAMU, College Station, TX 77843-3367, USA, jleon@tamu.edu}}
\date{}

\makeatletter
\newenvironment{breakablealgorithm}
{
	\begin{center}
		\refstepcounter{algorithm}
		\hrule height.8pt depth0pt \kern2pt
		\renewcommand{\caption}[2][\relax]{
			{\raggedright\textbf{\fname@algorithm~\thealgorithm} ##2\par}%
			\ifx\relax##1\relax 
			\addcontentsline{loa}{algorithm}{\protect\numberline{\thealgorithm}##2}%
			\else 
			\addcontentsline{loa}{algorithm}{\protect\numberline{\thealgorithm}##1}%
			\fi
			\kern2pt\hrule\kern2pt
		}
	}{
		\kern2pt\hrule\relax
	\end{center}
}
\makeatother

\newtheorem{definition}{Definition}
\newtheorem{theorem}{Theorem}
\newtheorem{lemma}{Lemma}

\newtheorem{proposition}{Proposition}

\theoremstyle{definition}
\allowdisplaybreaks

\algnewcommand\algorithmicforeach{\textbf{for each}}
\algdef{S}[FOR]{ForEach}[1]{\algorithmicforeach\ #1\ \algorithmicdo}

\begin{document}
	
	\newcolumntype{L}[1]{>{\raggedright\arraybackslash}p{#1}}
	\newcolumntype{C}[1]{>{\centering\arraybackslash}p{#1}}
	\newcolumntype{R}[1]{>{\raggedleft\arraybackslash}p{#1}}
	
	\maketitle
	\vspace{-1.5cm}
	\begin{abstract}	
		This paper studies binary linear programming problems in the presence of uncertainties that may cause solution values to change during implementation. This type of uncertainty, termed implementation uncertainty, is modeled explicitly affecting the decision variables rather than model parameters. The binary nature of the decision variables invalidates the use of the existing models for this type of uncertainty. The robust solutions obtained are optimal for a worst-case min-max objective and allow a controlled degree of infeasibility with respect to the associated deterministic problem. Structural properties are used to reformulate the problem as a mixed-integer linear binary program. The degree of solution conservatism is controlled by combining both constraint relaxation and cardinality-constrained parameters. Solutions for optimization problems under implementation uncertainty consist of a set of robust solutions; the selection of solutions from this possibly large set is formulated as an optimization problem over the robust set. Results from an experimental study in the context of the knapsack problem suggest the methodology yields solutions that perform well in terms of objective value and feasibility. Furthermore, the selection approach can identify robust solutions that possess desirable implementation characteristics.
		
	\end{abstract}
	
	\section{Introduction}\label{introsection}

\subsection{Introduction}

This paper studies binary linear programming problems in the presence of uncertainties that may cause solution values to change during implementation. This type of uncertainty, termed implementation uncertainty, is modeled explicitly affecting the decision variables rather than model parameters. Implementation uncertainty may result in implemented solutions that are different from what is prescribed by the BLP. The impact of implementation uncertainty on binary variables can be seen as if the variable is switching its prescribed value at the time of the implementation; therefore, a different solution than the prescribed one is implemented. Implementation uncertainty inevitably occurs due to inherent fidelity limitations of problem formulations and unexpected future events, including those caused by exogenous factors such as political directives, regulatory issues, or sudden extreme events. Model fidelity limitations are unavoidable in practice due to restricted time availability during modeling, limited knowledge about the problem at hand, and simplifying model assumptions. Implementing a different solution rather than the prescribed one may cause the objective value to become negatively impacted, leading to a suboptimal value, and the implemented solution may no longer be feasible. This type of uncertainty affecting binary variables hinders applying most of the existing uncertainty models proposed in the related literature. Assuming that only a subset of variables is affected by implementation uncertainty and the others are deterministic, solving a BLP under implementation uncertainty can be viewed as specifying the values on the deterministic variables while the uncertain ones may take any possible value. Consider the following simple example to illustrate the proposed problem addressed in this paper.

\subsection{Illustrative Example}

An investor needs to decide which of ten projects to invest in. Each project $i$ has associated profit $c_i$ and cost, $a_i$. The selection has to be such that it maximizes the profit while maintaining the cost within a budget $b=26$; the remaining data for this example is given in Table \ref{tabexample}:

\begin{table}[h]
	\centering
	\begin{tabular}{| c | c  c  c  c  c  c  c  c  c  c |}
		\hline
		$i$ & 1 & 2 & 3 & 4 & 5 & 6 & 7 & 8 & 9 & 10 \\
		\hline
		$c_i$ & 7& 3 & 9 & 9 & 10 & 7 & 4 & 2 & 6 & 2 \\
		$a_i$ & 4 & 5 & 9 & 8 & 4 & 4 & 6 & 6 & 2 & 3 \\ \hline
	\end{tabular}
	\caption{Profits and costs for the candidates projects in the illustrative example.}
	\label{tabexample}
\end{table}

The optimal solution for this deterministic version of the problem is a profit $Z_{Det^*}=41$ and a total budget requirement of $LHS=26=b$.

To introduce the concept of implementation uncertainty, consider that projects 1 and 2 are from a cloud-based investment market where investment opportunities are known to be unpredictable as projects may become unavailable between the time the investment decision is made and the time the decision is implemented. Management is interested in considering projects from this cloud market (i.e., projects 1 and 2) in addition to their regular investment opportunities (i.e., projects 3 to 10). In this paper, we say that implementation uncertainties affect the projects in the cloud market and that the associated implemented decision may change from the decisions described by the optimization model; we call these variables "uncertain"; furthermore, we assume that variables associated with regular projects will not change from the definitions described by the optimization model; we call these variables "deterministic." Because of implementation uncertainty, the objective value and the feasibility of a prescribed solution can be affected. Let $x_i$ be the binary decision variable associated with the selection of project $i$, Table \ref{taboutcomes} shows all possible outcomes if the deterministic optimal solution were to be implemented as prescribed; i.e., variables 3-10 values would remain as specified, while variables 1 and 2 may take any value during implementation. Hence, in this paper, a solution to a BLP under implementation uncertainty is a set of solutions specified by the values of the deterministic variables. In this paper, we answer the question: Is there a robust solution set with the desired objective values that can guarantee the desired feasibility level? As will be described in Section \ref{solselection} of the paper, once this set is determined, one may formulate an optimization problem over the robust solution set to single out a solution for implementation.


\begin{table}[h]
	\centering
	\begin{tabular}{| c | c  c | c  c  c  c  c  c  c  c | c  c |}
		\hline
		$i$ & 1 & 2 & 3 & 4 & 5 & 6 & 7 & 8 & 9 & 10 & $Z$ & LHS \\
		\hline
		\multirow{4}{*}{$x_i$} & 0 & 0 & & & & & & & & & 34 & 22  \\
		& 0 & 1 & & & & & & & & & 37 & \textit{27} \\ 
		& \textbf{1} & \textbf{0} & \textbf{1} & \textbf{0} & \textbf{1} & \textbf{1} & \textbf{0} & \textbf{0} & \textbf{1} & \textbf{1} & \textbf{41} & \textbf{26} \\
		& 1 & 1 & & & & & & & & & 44 & \textit{31} \\ \hline
	\end{tabular}
	\caption{All possible outcomes after implementing the optimal deterministic solution. The prescribed solution is shown in bold font; italic font indicates infeasible with respect to the deterministic.}
	\label{taboutcomes}
\end{table}

In Table \ref{taboutcomes}, the worst-case profit occurs in the outcome where both projects 1 and 2 are canceled, yielding a lower profit of 34 instead of the prescribed 41. In terms of feasibility, two of the four possible outcomes can be infeasible due to uncertainty. If projects 1 and 2 had been part of the initial optimal solution, any unexpected change in their values would result in profit loss, and the outcome would always be feasible. On the other extreme, if neither project were in the deterministic solution, then any unexpected change in their values would cause an increase in the profit and a possible violation of feasibility. In general, implementation uncertainty may cause infeasibilities during implementation; hence, when devising solution methodologies to tackle optimization problems under this type of uncertainty, the solution approach must explicitly consider possible feasibility violations (with respect to the deterministic problem) – in this paper, the level of acceptable infeasibility is controlled by adding a feasibility parameter $\delta_j \geq 0$ in the RHS of every constraint.

The problem studied in this paper aims at finding solutions that have close-to-optimal profit and can guarantee a desirable level of feasibility. Table \ref{tabrobust} displays the solution set found using the methodology in Section \ref{MIPRob} allowing a maximum of 5\% violation of the budget constraint, and a method in Section \ref{exps3} was used to specify a specific solution ($\hat{x}_{\delta}^D$) from the robust set. Compared to the deterministic solution, the robust solution displays lower profit than the deterministic optimal value (39 versus 41) but better feasibility performance (i.e., 3 out of 4 versus 2 out of 4 feasible outcomes). One can note that the proposed methodology was able to find a feasible solution where project 1 switches from 1 to 0 and project 2 is selected (i.e., switching from 0 to 1) - this situation was infeasible when using the optimal deterministic solution (see Table \ref{taboutcomes}).

\begin{table}[h]
	\centering
	\begin{tabular}{| c | c  c | c  c  c  c  c  c  c  c | c  c |}
		\hline
		$i$ & 1 & 2 & 3 & 4 & 5 & 6 & 7 & 8 & 9 & 10 & $Z$ & LHS \\
		\hline
		\multirow{4}{*}{$x_i$} & 0 & 0 & & & & & & & & & 32 & 18  \\
		& 0 & 1 & & & & & & & & & 35 & 23 \\ 
		& \textbf{1} & \textbf{0} & \textbf{0} & \textbf{1} & \textbf{1} & \textbf{1} & \textbf{0} & \textbf{0} & \textbf{1} & \textbf{0} & \textbf{39} & \textbf{22} \\
		& 1 & 1 & & & & & & & & & 42 & \textit{27} \\ \hline
	\end{tabular}
	\caption{Robust solution obtained with the proposed method. The prescribed solution is shown in bold font; italic font indicates infeasible with respect to the deterministic.}
	\label{tabrobust}
\end{table}

Of notice is that the feasibility of the robust solution seems more “controlled,” as evidenced by a maximum feasibility violation of 1 unit versus 5 units for the deterministic solution. Arguably, in many cases, a “small” feasibility violation would be “easier” to fix during implementation. In our example, if the infeasibility were to occur during implementation, it would be easier for the company to request a small increase in the budget to accommodate a good investment opportunity, e.g., request an increase in one unit of budget achieving a profit of 42. Thus, although not always true, in this example, the company could achieve better profit than the deterministic optimal by fixing the small infeasibility that occurred at implementation.


Another feature of the proposed approach is that, if they exist, it can generate solutions that would guarantee feasibility during implementation. Although, these worst-case solutions may come at significant degradation of the objective value compared to the deterministic optimal. Situations in practice where these extreme worst-case solutions are required include decisions that relate to safety, loss of life, or negative impact to large populations of people; e.g., electricity network reliability, design of vaccine distribution systems, complex surgery procedure planning, design of anti-terrorism systems, warfare planning, etc.

\subsection{Literature Review}

Different approaches aim to protect the optimality and feasibility of solutions in the face of uncertainties including stochastic optimization \citep[e.g.][]{dantzig1955linear, beale1955minimizing, wets1966programming, wets1974stochastic, wets1983solving} and robust optimization \citep[e.g.][]{soyster1973technical, mulvey1995robust,bertsimas2004price}. Stochastic optimization seeks solutions that remain optimal and feasible with high probability. However, there may exist realizations of the uncertainty where the optimality or feasibility are not satisfied \citep[see][]{ben2009robust}. On the other hand, robust optimization approaches seek solutions that satisfy the given levels of optimality and feasibility for any realization of the uncertainty; such solutions are termed \textit{robust solutions} \citep[][]{mulvey1995robust}. For instance, \cite{soyster1973technical} considers perturbations in the coefficients of the constraints using convex sets; the resulting model produces solutions that are feasible for any realization of the data within the convex sets.

The existing work in the field of robust optimization accounting for implementation uncertainty is very limited \citep[][]{gabrel2014recent}. \cite{ben2009robust} propose two forms of modeling implementation uncertainty on real decision variables: additive implementation errors refer to the case when a random value is added to the prescribed value, and multiplicative implementation error refers to the case when the random value multiplies the prescribed value; furthermore, the authors show that these forms of implementation errors are equivalent to artificial data uncertainties and can be treated as such. These forms of modeling implementation uncertainty have been used in single optimization problems \citep[e.g.][]{das1997nonlinear,lewis2009lipschitz},  and in multiobjective optimization problems \citep[e.g.][]{deb2006introducing,jornada2016biobjective,eichfelder2017decision}. However, these models of implementation uncertainty cannot be extended to the case of binary problems when dealing with implementation uncertainty because using equivalent data uncertainty is not straightforward; additional discussion is provided in Section \ref{modelsection}.


The methodologies in this paper assume worst-case objective functions \citep[][]{kouvelis1997robust}; the study of models with other types of objective functions remains as future research. A characteristic of worst-case robust solutions is that they tend to be conservative because they may excessively sacrifice optimality to satisfy the given level of feasibility. For instance, the model in \cite{soyster1973technical} is considered too conservative from this perspective \citep[][]{bertsimas2004price}. Different authors have addressed this issue by modeling uncertainty using different representations; for instance, \cite{el1997robust,el1998robust,ben1998robust,ben1999robust,ben2002robustmethod} propose a less conservative model by using ellipsoidal sets to describe data uncertainty, \cite{bertsimas2003robust, bertsimas2004price} control conservatism by bounding the maximum number of uncertain coefficients changing in each constraint simultaneously, and \citep[][]{kouvelis1997robust} use measures of robustness that seek to minimize the difference between the objective and robust objective values (i.e., the maximum deviation).

To address the conservatism of the robust solutions, the proposed robust formulation includes a feasibility parameter that gives the decision-maker direct control of the feasibility level to improve the robust solutions' objective performance. In addition, the proposed methodology also incorporates cardinality-constrained concepts inspired by the work in \cite{bertsimas2004price} applied to implementation uncertainty directly affecting decision variables rather than model parameters.

The rest of the paper is organized as follows. Section \ref{modelsection} describes the robust optimization under implementation uncertainty model; Section \ref{robust} presents the formulation of the problem and characteristics of robust solution sets; Section \ref{solution} describes the solution methodology to find robust solution sets and how to select solutions for implementation; Section \ref{expsection} presents an experimental study in the context of the Knapsack Problem; and Section \ref{conclu} contains concluding remarks.
	\section{Model Development}\label{modelsection}

Let $\mathbb{B}^n \subset \mathbb{R}^n$ be the set of $n$-dimensional binary vectors $\mathbb{B}^n = \left\{x = (x_1, ..., x_n) : \right.$ $\left. x_i \in \{0,1\}, \forall i \in I \right\}$, where $I=\{1,...,n\}$. Let $x\in \mathbb{B}^n$ be the vector of decision variables $x_i$ and let $f:\mathbb{B}^n \rightarrow \mathbb{R}$ be a function defined as $f(x) = \sum_{i=1}^{n}c_ix_i$, with $c_i \in \mathbb{R}, \forall i \in I$. Let $g_j:\mathbb{B}^n \rightarrow \mathbb{R}$ be a function defined as $g_j(x) = \sum_{i=1}^{n}a_{ij}x_i$ with $a_{ij} \in \mathbb{R}, \forall i \in I,j \in J=\{1,...,m\}$ defining the left-hand-side, and $b_j \in \mathbb{R}$ defining the right-hand-side of the $j$-th constraint $\forall j \in J$. The feasible set $X$ is defined as $X=\left\{x \in \mathbb{B}^n: g_j(x) \leq b_j, \forall j\in J\right\}$. A binary linear programming problem (BLP) can be formulated as follows:
\begin{equation}\label{blp}
	\min_{x \in \mathbb{B}^n} \{f(x): x \in X\}.
\end{equation}

Henceforth, formulation \eqref{blp} is termed \textit{the deterministic BLP formulation}, and $X$ is termed \textit{the deterministic feasible set}.

Let $\hat{x}$ denote a prescribed solution and $\tilde{x}$ denoted an implemented solution. It is assumed that $\hat{x}$ is the solution obtained from solving an optimization model, and $\tilde{x}$ is the actually implemented solution, which may be different than $\hat{x}$ due to implementation uncertainty.

\begin{definition}\label{implementationerror}
	A binary variable $x_i$ is under implementation uncertainty if $p_i < 1$ or $q_i < 1$, and the following conditional probabilities hold true:	
	\begin{equation}\label{probuncertain}
		\begin{split}
			P(\tilde{x}_i = \hat{x}_i | \hat{x}_i = 0) = p_i, \quad\quad & P(\tilde{x}_i = 1 - \hat{x}_i | \hat{x}_i = 0) = 1 - p_i, \\
			P(\tilde{x}_i = \hat{x}_i | \hat{x}_i = 1) = q_i, \quad\quad & P(\tilde{x}_i = 1 - \hat{x}_i | \hat{x}_i = 1) = 1 - q_i.
		\end{split}		
	\end{equation}	
\end{definition}

\cite{ben2009robust} propose the additive and multiplicative implementation errors to model implementation uncertainty in real variables. The additive implementation errors consists of a random value $\epsilon$ added to the prescribed value of the decision variable, $\tilde{x} = \hat{x}+\epsilon$. The multiplicative implementation errors a random value multiplies the prescribed value of the decision variable, $\tilde{x}=\epsilon \hat{x}$. In the two models, the value of $\epsilon$ belongs to a defined set. These models of implementation uncertainty cannot handle the case of binary variables appropriately. For instance, in the case of additive implementation errors, and assuming $epsilon \in \left\{-1,0,1\right\}$, adding $\epsilon$ to the value decision variable may generate infeasible values for the decision variable. For example, with $\hat{x} = 0$ the implemented value $\tilde{x}$ may be -1, 0 or 1; similarly, with $\hat{x} = 1$ the implemented value $\tilde{x}$ may be 0, 1 or 2. On the other hand, in the case of multiplicative implementation error with $\hat{x} = 0$ the result is $\tilde{x} = \epsilon \hat{x} = 0$ for any value of $\epsilon$; hence, these model of implementation uncertainty is not appropriate neither.

Lemma \ref{lineareq2} in Section \ref{MIPRob} presents a model to handle the impact of implementation uncertainty in binary variables appropriately. The rest of this section presents several concepts, assumptions, and notation used throughout the paper.

Variables affected by implementation uncertainty are termed \textit{uncertain variables}, otherwise they are termed \textit{deterministic variables}. Without loss of generality, the decision vector $x$ is decomposed into two vectors $x_C$ and $x_U$, where $x_C$ is composed of the deterministic variables $x_1,..., x_c$, and $x_U$ is composed of the uncertain variables $x_{c+1},..., x_n$; for convenience, define $C = \{1,...,c\}$ as the set of indices of the deterministic variables in $x_C$, and $U =\{c+1,...,n\}$ as the set of indices of the uncertain variables in $x_U$.


\noindent \textbf{Assumption 1.} Sets $C$ and $U$ are given; i.e., it is known which variables are deterministic and uncertain.

\noindent \textbf{Assumption 2.} The probabilities $p_i$ and $q_i$ for the uncertain variables are unknown; i.e., it is known that uncertain variables may have a different implemented value, but the probability that a change may occur is unknown.



For a prescribed solution $\hat{x} = (\hat{x}_C, \hat{x}_U)$, a corresponding implemented solution $\tilde{x} = (\tilde{x}_C, \tilde{x}_U)$ has the same values of the deterministic variables, $\tilde{x}_C = \hat{x}_C$, and the value $\tilde{x}_U$ possibly taking any combination of the $n-c = u$ uncertain variables. The \textit{set of implemented outcomes}, $\mathcal{U}(x)$, is defined as follows:

\begin{definition}\label{uncertainsetbinary}
	Given $x \in \mathbb{B}^n$, the set of implemented outcomes associated with $x$ is defined as $\mathcal{U}(x) = \left\{\tilde{x}=(\tilde{x}_C,\tilde{x}_U)\in \mathbb{B}^n : \tilde{x}_C = x_C \right\}$.
\end{definition}

The number of possible implementation outcomes grow exponentially with the number of uncertain variables; i.e., $|\mathcal{U}(x)| = 2^{|U|}$. Furthermore, $\exists \tilde{x} \in \mathcal{U}(x)$ such that $\tilde{x} = x$.

	\section{RBIU-$\delta$ Problem Formulation}\label{robust}

When impacted by implementation uncertainty, $\tilde{x}$ may result in objective function values different than the objective function value of $\hat{x}$ or even become infeasible. A robust BLP under implementation uncertainty (RBIU-$\delta$) aims at finding solutions that guarantee desired levels of optimality and feasibility in the face of implementation uncertainty.

The objective robustness level, $\gamma(x)$, measures the degree to which a solution's objective function value degrades when affected by implementation uncertainty.

\begin{definition}\label{absoluterobustness}
	Given a binary vector $x$, the objective robustness level, $\gamma(x)$, is defined as:	
	\begin{equation}\label{absform}
		\gamma(x) = \max_{y \in \mathcal{U}(x)} \left \{f(y) \right \}.
	\end{equation}
\end{definition}

$\gamma(x)$ provides the worst-case value of the objective function among all outcomes in $\mathcal{U}(x)$, guaranteeing that the objective value of the implemented solution does not worsen when affected by implementation uncertainty. We use the worst-case objective for convenience as it allows us to linearize \eqref{absform}.

\begin{lemma}\label{multiplesolutions}
	Given $x = (x_C,x_U) \in \mathbb{B}^n$, then $\gamma(y) = \gamma(x), \forall y \in \mathcal{U}(x)$.
\end{lemma}

\begin{proof}
	By Definition \ref{uncertainsetbinary} it follows that $y_C = x_C$, and therefore $\mathcal{U}(y) = \mathcal{U}(x)$. From Definition \ref{absoluterobustness}, every vector $y \in \mathcal{U}(x)$ produces the same objective robustness value $\gamma(x)$.
\end{proof}


Typical in robust optimization is that solutions tend to be very conservative (i.e., the resulting degradation in objective function may be too excessive), in particular, if one desires to preserve feasibility with respect to the deterministic version of the problem (we term this type of feasibility \textit{deterministic-feasibility}).  One can obtain a better objective performance in robust optimization by accepting some degree of infeasible outcomes. For instance, when using cardinality-constrained approaches \citep[see][]{bertsimas2004price}, $CC$, infeasible outcomes occur when more than $\Gamma$ coefficients are affected by uncertainty – Successful application of $CC$ approaches rely on proving a low probability of the occurrence of infeasibilities. Unfortunately, $CC$ does not explicitly offer any guarantee in the level of deterministic feasibility violation, which may be problematic when dealing with binary variables affected by implementation uncertainty. To avoid this problem, similar to \cite{ben2000robust} concept of “feasibility tolerance,” we introduce the feasibility parameter, $\delta_j \geq 0$, to guarantee a maximum level of constraint violation while simultaneously expanding the robust solution space to include solutions with better objective performance. 
The deterministic constraints are reformulated as follows:

\begin{equation}\label{feaslevelcons}
	\max_{y \in \mathcal{U}(x)}\left\{g_j(y)\right\} \leq b_j + \delta_j, \quad \forall j \in J.
\end{equation}

The RBIU-$\delta$ is formulated as follows:
\begin{equation}\label{RBIU}
	\min_{x\in \mathbb{B}^n} \left \{\gamma(x):  x \in \mathcal{X}\right \};
\end{equation}

\noindent where $\mathcal{X} \subset \mathbb{B}^n$ denotes the feasible region defined by the modified constraints \eqref{feaslevelcons}. $\mathcal{X} \subset X$, $\mathcal{X} \cap X = \emptyset$ and $\mathcal{X} = \emptyset$ are possible; furthermore, the intersection may exist and not contain the optimal deterministic solution.

\begin{figure}[h]
	\centering \includegraphics[scale = 0.6]{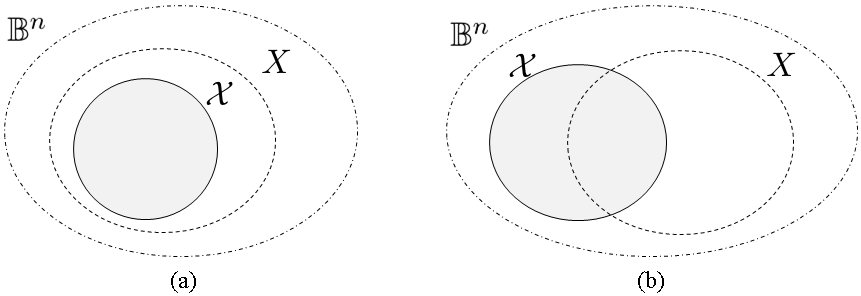}
	\caption{Illustration of the feasible sets $X$ and $\mathcal{X}$. (a) Illustrates the most conservative case when $\delta_j = 0, \forall j \in J$. (b) Illustrates a less conservative case when $\exists j$ such that $\delta_j >0$. \label{plot06}}{}
\end{figure}

Setting $\delta_j = 0, \forall j \in J$ guarantees robust solutions, if they exist, that are also feasible with respect to the deterministic problem; i.e., $\mathcal{X} \subseteq X$ (Figure \ref{plot06}(a)). However, this feasibility guarantee may come at the expense of objective function value degradation as the deterministic optimal solution may not belong to $\mathcal{X}$ leading to a robust solution with a worse objective function value. The case for $\delta_j = 0, \forall j \in J$ is akin to conservative robust approaches such as the one proposed by \cite{soyster1973technical}. Setting  $\delta_j>0$ expands the robust feasible region $\mathcal{X}$ to include solutions that may or may not be deterministic-feasible after implementation. As it becomes clear in Figure \ref{plot06}(b), a distinctive property of the proposed $\delta-$approach is that it will evaluate solutions that might reside outside of $X$ that may have a more desirable outcome when affected by implementation uncertainty when compared to CC that is restricted to solutions in $X$. In summary, $\Gamma$ and $\delta$ have a similar effect on achieving improved objective performance at the expense of bounded feasibility violations; additionally, as it is more critical in binary formulations, $\delta$ offers a guarantee on the maximum level of deterministic-feasibility violation. The proposed RBIU-$\delta$-CC \eqref{cbiu} takes advantage of both methods, giving the decision maker the ability to control the degree of conservativism of the solutions while at the same time achieving further improvements in objective value performance while guaranteeing a level of deterministic constraint violation. The gains obtained by combining both approaches are evident in the experimental results.


\begin{lemma}\label{allsolutions}
	Let $x^* = (x^*_C,x^*_U)$ be an optimal solution to the RBIU-$\delta$. Any solution $y = (y_C,y_U) \in \mathcal{U}(x^*)$ is also optimal to the RBIU-$\delta$.
\end{lemma}

\begin{proof}
	By Lemma \ref{multiplesolutions}, $\gamma(y) = \gamma(x^*)$. $\mathcal{U}(x^*) = \mathcal{U}(y)$ because $x^*_C = y_C$, and $x^*_U$ and $y_U$ can take any combination of the uncertain variables using Definition \ref{uncertainsetbinary}; therefore, $y \in \mathcal{X}$. Hence, $y$ is also optimal for the RBIU-$\delta$.
\end{proof}

As implied by Lemma \ref{allsolutions}, a characteristic of the RBIU-$\delta$ is that it has multiple optimal solutions (all solutions in $\mathcal{U}(x^*)$), and that to find the set of solutions $\mathcal{U}(x^*)$ it is only necessary to find one solution in the set suggesting a solution methodology that only determines $x^*_C$.

The set of solutions of the RBIU-$\delta$, $\mathcal{U}(x^*)$, is named here the \textit{robust-optimal solution set}, $\mathcal{U}^*$. If $\mathcal{U}^*$ exists, contains a solution for each combination of the uncertain variables.
	\section{Solution Methodology}\label{solution}

The methodology to solve \eqref{RBIU} consists of two stages. Stage I aims at finding a robust-optimal solution set, $\mathcal{U}^*$; given the possibly large size of the optimal set, (Sections \ref{MIPRob} and \ref{cardinality}). Stage II aims at selecting desirable properties from $\mathcal{U}^*$, (Section \ref{solselection}).

\subsection{Mixed-Integer Linear Programming Reformulation, RBIU-$\delta$}\label{MIPRob}

When solving \eqref{RBIU}, the evaluation of \eqref{absoluterobustness} and \eqref{feaslevelcons} can be burdensome because they require finding maximum values among all vectors in $\mathcal{U}(x)$ for every $x \in \mathcal{X}$. Lemma \ref{lineareq2} provides equivalent linear formulations of \eqref{absoluterobustness} and \eqref{feaslevelcons} allowing to find a Stage I optimal solution of the RBIU-$\delta$ by solving a mixed-integer linear programming problem (MILP). Advantages of the MILP reformulation are: 1) a linearization of the RBIU-$\delta$, and 2) a reduction of the search space from $2^n$ to $2^{|C|}$, where $|C| < n$.

\begin{lemma}\label{lineareq2}
	Expression $\sum_{i \in C}c_{i}x_i + \sum_{i \in U}(c_{i} + |c_{i}|)/2$ and $\sum_{i \in C}a_{ij}x_i + \sum_{i \in U}(a_{ij} + |a_{ij}|)/2$ are equivalent to $\max_{y \in \mathcal{U}(x)} \left \{f(y) \right\}$ and $\max_{y \in \mathcal{U}(x)} \left \{g_j(y) \right\}, \forall j \in J$, respectively.
\end{lemma}

\begin{proof}
	Expression $(c_{i} + |c_{i}|)/2 > c_ix_i$ for any value of $x_i$. Note that if $c_{i} \geq 0$, then $|c_{i}| = c_{i}$ and $(c_{i} + |c_{i}|)/2 = (c_{i} + c_{i})/2= c_{i}\geq c_{i}x_i$; similarly, if $c_{i} < 0$, then $|c_{i}| = -c_{i}$ and $(c_{i} + |c_{i}|)/2 = (c_{i} - c_{i})/2 = 0 \geq c_{i}x_i$. Therefore, $(c_{i} + |c_{i}|)/2 \geq c_{i}x_i$ for any value of $x_i$. Expression $\max_{y \in \mathcal{U}(x)} \left \{f(y)\right\}$ is rewritten as $\max_{y \in \mathcal{U}(x)} \left \{\sum_{i \in C}c_{i}y_i + \sum_{i \in U}c_{i}y_i\right\} = \max_{y \in \mathcal{U}(x)} \left \{\sum_{i \in C}c_{i}y_i\right\} + \max_{y \in \mathcal{U}(x)} \left \{\sum_{i \in U}c_{i}y_i\right\}$. Given that values of the deterministic variables $x_i, i \in C$ are fixed, then $\max_{y \in \mathcal{U}(x)} \left \{\sum_{i \in C}c_{i}y_i\right\} = \sum_{i \in C}c_{i}x_i$. On the other hand, the value of $\max_{y \in \mathcal{U}(x)} \left \{\sum_{i \in U}c_{i}y_i\right\}$ depends of the value $y_i, i \in U$, and from the previous result it follows that $\max_{y \in \mathcal{U}(x)} \left \{\sum_{i \in U}c_{i}y_i\right\} = \sum_{i \in U}\max_{y \in \mathcal{U}(x)} \left \{c_{i}y_i\right\} = \sum_{i \in U}(c_{i} + |c_{i}|)/2$ for any value of $y_i, i \in U$. Therefore, $\max_{y \in \mathcal{U}(x)} \left \{f(y) \right\} = \sum_{i \in C}c_{i}x_i + \sum_{i \in U}(c_{i} + |c_{i}|)/2$. It can be proved similarly that $\max_{y \in \mathcal{U}(x)} \left \{g_j(y) \right\} = \sum_{i \in C}a_{ij}x_i + \sum_{i \in U}(a_{ij} + |a_{ij}|)/2, \forall j \in J$.
\end{proof}

The RBIU-$\delta$ is a MILP reformulation of \eqref{RBIU} as follows:
\begin{equation}\label{rmbp}
	\min_{x \in \mathbb{B}^c}\left\{\gamma(x): x \in \mathcal{X}^{'} \right\};
\end{equation}

\noindent where $\mathcal{X}^{'}$ is the feasible region such that expressions $\max_{y \in \mathcal{U}(x)} \left \{f(y) \right \}$ and $\max_{y \in \mathcal{U}(x)} \left \{g_j(y) \right\}, \forall j \in J$ are replaced by their equivalent expressions in Lemma \ref{lineareq2}. The linear reformulations depend on deterministic variables only because the uncertain variables are replaced by constant values; therefore, the search space of the RBIU-$\delta$ is $\mathbb{B}^{|C|}$, with $|C| < n$.

Formulation \eqref{rmbp} is used to find, if it exists, one solution of the RBIU-$\delta$, $x^* = (x^*_C,x^*_U)$, and therefore, the robust-optimal solution set $\mathcal{U}^*$.

\subsection{Combined Feasibility Parameter and Cardinality-Constrained Robust Formulation, RBIU-$\delta$-CC}\label{cardinality}

A second approach provides additional control of the solution's conservatism by controlling the maximum number of variables, $\Gamma$, assumed to be affected by uncertainty when solving the robust formulation in addition to the feasibility parameter $\delta_j$. We term this approach \emph{combined feasibility parameter and cardinality constrained robust formulation} for a BLP under implementation uncertainty (RBIU-$\delta$-CC). In contrast to Bertsimas and Sim's robust formulation, RBIU-$\delta$-CC considers uncertainty affecting the variables instead of the coefficients of the model and impacting the entire column of the uncertain variables, including constraints and the objective function simultaneously. Moreover, the interval of uncertainty in the proposed model is asymmetric and binary due to the nature of the decision variables. The RBIU-$\delta$-CC is based on the linearization model in Lemma \ref{lineareq2}.

\subsubsection{RBIU-$\delta$-CC Formulation}\label{devcc}

Consider the expressions in Lemma \ref{lineareq2} for maximum contribution of $a_{ij}x_i$ and $c_ix_i$. Let $\Gamma$ be an integer control parameter name de \textit{cardinality-constrained parameter} such that with $1 \leq \Gamma \leq |U|$; $\Gamma$ represents the maximum number of uncertain variables that may have different prescribed and implemented values. The RBIU-$\delta$-CC is formulated as follows:
\begin{equation}\label{cbiu}
	\min_{x \in \mathbb{B}^n}\left\{\gamma(x): x \in \mathcal{X}^{''} \right\}.
\end{equation}

Where $\mathcal{X}^{''}$ is a feasible region defined as follows:
\begin{alignat}{3}
	& \sum_{i \in C}c_i x_i + \max_{\left \{S_0 : S_0 \subseteq U, |S_0| \leq \Gamma \right \}}\left\{ \sum_{i \in S_0} \left(\frac{c_i + |c_i|}{2}\right) + \sum_{i \in U \setminus S_0 }c_i x_i \right\} \leq \gamma(x), \label{maxsolcard} \\
	& \sum_{i \in C}a_{ij} x_i + \max_{\left \{S_j : S_j \subseteq U, |S_j| \leq \Gamma \right \}}\left\{\sum_{i \in S_j} \left (\frac{a_{ij} + |a_{ij}|}{2} \right ) + \sum_{i \in U \setminus S_j}a_{ij} x_i \right\}  \leq b_j + \delta_j, \forall j \in J. \label{maxmodcard}
\end{alignat}

RBIU-$\delta$-CC seeks for the combination of at most $\Gamma$ uncertain variables affected by implementation uncertainty producing the maximum degradation of the objective function value and the maximum value of the constraints' uncertain component that satisfies the desired feasibility level, and protects the optimality and feasibility levels against any, at most, $\Gamma$ uncertain variables with different prescribed and implemented values. 

The selection of the at most $\Gamma$ uncertain variables requires the enumeration of all the subsets $S_0$ and $S_j$ of $U$ whose cardinality is less than or equal to $\Gamma$. The development of an equivalent linear reformulation to the RBIU-$\delta$-CC follows the work in \cite{bertsimas2004price}. The RBIU-$\delta$-CC is equivalent to the following linear reformulation:

\begin{align}\label{cmbp}
	\min \quad & \gamma(x) \\
	\text{s.t. } \quad & \sum_{i \in C}c_i x_i + \Gamma v_0 + u_{00} + \sum_{i \in U}u_{i0} \leq \gamma \\
	& v_0 + u_{i0} \geq \frac{c_i + |c_i|}{2} - c_ix_i, & \quad \forall i \in U \\
	& u_{00} \geq \sum_{i \in U}c_ix_i \\
	& \sum_{i \in C}a_{ij} x_i + \Gamma v_j + u_{0j} + \sum_{i \in U}u_{ij} \leq b_j + \delta_j, & \quad \forall j \in J\\
	& v_j + u_{ij} \geq \frac{a_{ij} + |a_{ij}|}{2} - a_{ij}x_i, & \quad \forall i \in U, j  \in J \\
	& u_{0j} \geq \sum_{i \in U}a_{ij}x_i, & \quad \forall j \in J \\
	& u_{i0}, u_{ij}, v_0, v_j \geq 0, & \quad \forall i \in U, j \in J \\
	& u_{00}, u_{0j} \text{ are unrestricted} & \quad \forall j \in J \\
	& x_i \in \{0, 1 \}, & \forall i \in I.
\end{align}

The equivalent linear reformulation of the RBIU-$\delta$-CC contains: 1) $n$ binary variables, 2) $m+1$ unrestricted variables, 3) $m|U| + m + |U| + 1$ nonnegative variables, and 4) $m|U|+2m+|U|+n+3$ constraints.

RBIU-$\delta$-CC produces more optimistic solutions than RBIU-$\delta$ at the expense of possible infeasibilities to the desired feasibility level or degradation of the optimality. The following section presents an upper bound of the probability of infeasibilities or optimality degradation of the RBIU-$\delta$-CC solutions.

Of notice is that RBIU-$\delta$-CC is general in the sense that it becomes RBIU-$\delta$ if $\Gamma=|U|$, and works as a cardinality-constrained approach for $\delta = 0$. Additionally, the result $\mathcal{U}(x) = \mathcal{U}(y), \forall y \in \mathcal{U}(x)$ is not true for the RBIU-$\delta$-CC formulation because in this formulation only a subset of uncertain variables is being impacted by implementation uncertainty simultaneously, except when $\Gamma = |U|$.

\subsubsection{Probability Bounds}\label{bounds}

RBIU-$\delta$-CC may become infeasible if more than $\Gamma$ uncertain variables change their values during implementation; this section presents derivations of upper bounds on the probability that this type of infeasibility occurs. To estimate these probability upper bounds, Assumption 2 is relaxed  and we assume that $p_i = p$ and $q_i = q$ ,$\forall i \in I$, with $p$ and $q$ known.

Let $\eta_0$ and $\eta_1$ be two independent random variables such that $\eta_0$ measures the number of uncertain variables with prescribed value 0 and implemented value 1, and $\eta_1$ measures the number of uncertain variables with prescribed value 1 and the implemented value 0. Then, the probability that there exist exactly $\Gamma$ uncertain variables with different prescribed and implemented values can be computed as follows:
\begin{equation}\label{distofchanging}
	P(\eta_0 + \eta_1 = \Gamma) = \sum_{i = 0}^{\Gamma}\binom{U_0}{i}(1-p)^ip^{U_0 - i}\binom{U_1}{\Gamma-i}(1-q)^{\Gamma-i}q^{U_1 - \Gamma + i};
\end{equation}

\noindent where $U_0$ and $U_1$ are the number of uncertain variables whose prescribed values are 0 and 1, respectively, and $U_0 + U_1 = |U|$. Distribution \eqref{distofchanging} represents the sum of two independent binomial random variables $\eta_0$ and $\eta_1$.

\begin{theorem}\label{lowerprob}
	Given an RBIU-$\delta$-CC solution $x^*$ with optimality robustness level $\gamma(x^*)$, the upper bound of the probability that the objective value is not $\gamma(x^*)$ or any of the $j$-th model robustness constraints is not satisfied is given by:	
	\begin{equation}\label{ubprob}
		\begin{split}
			P\left(\left\{\sum_{i=1}^{n}c_ix_i^* > \gamma(x^*) \right\} \cup \left\{\sum_{i=1}^{n}a_{ij}x_i^* > b_j + \delta_j \right\}\right) \quad \quad\quad\quad\quad\quad \quad\quad\quad \quad\quad\quad \\ \leq 1 - \sum_{\ell = 0}^{\Gamma}\left(\sum_{i = 0}^{\ell}\binom{U_0}{i}(1-p)^ip^{U_0 - i}\binom{U_1}{\ell-i}(1-q)^{\ell-i}q^{U_1 - \ell + i} \right) .
		\end{split}
	\end{equation}	
\end{theorem}

\begin{proof}
	Given that the objective function value and constraints are protected for at most $\Gamma$ uncertain variables with different prescribed and implemented values, one can assume that they may become infeasible if there exists at least one more uncertain variable impacted by implementation uncertainty; this is $\eta_0 + \eta_1 > \Gamma$ in expression \eqref{distofchanging}. Therefore, expression \eqref{ubprob} holds.
\end{proof}

\subsection{Selected Robust 	Optimal Solutions}\label{solselection}

In Stage II, the problem of selecting a prescribed solution $\hat{x}$ for implementation is formulated as an optimization problem over the robust-optimal solution set. Given $\mathcal{U}^*$ generated by a solution $x^*=(x_C^*,x_U^*)$, the selection problem assumes the desired characteristics of the solution to be selected and recasts the deterministic version of the problem such that the deterministic variables are set to $x_C^*$ to ensure membership in the robust-optimal solution set $\mathcal{U}^*$. Notice that $\mathcal{U}^*$ can be obtained by solving RBIU-$\delta$ or RBIU-$\delta$-CC. As a result, the selected solution will have the desired robustness in terms of levels of objective value and feasibility performance and other characteristics that make it desirable for implementation. Examples of prescribed solutions $\hat{x} \in \mathcal{U}^*$ are presented in Table \ref{tab01}.

\begin{table}[h]
	\centering
	\begin{tabular}{| c | p{8cm} |}
		\hline
		Prescribed solution $\hat{x}$ & Selection Problem, SP \\ \hline
		$\hat{x}^D$ & $\text{SP}_1$: Solve \eqref{blp}; set $x_C \leftarrow x_C^*$ \\
		$\hat{x}^R$ & $\text{SP}_2$: Solve \eqref{blp}; set $x_C \leftarrow x_C^*$  and the right-hand side $\leftarrow b + \delta_j, \forall j \in J$\\
		$\hat{x}^{UB}$ & $\text{SP}_3$: Find $x\in \mathcal{U}^*:f(x) \geq f(y), \forall y \in \mathcal{U}^*$ \\
		$\hat{x}^{LB}$ & $\text{SP}_4$: Find $x\in \mathcal{U}^*:f(x) \leq f(y), \forall y \in \mathcal{U}^*$ \\ \hline
	\end{tabular}
	\caption{Selected prescribed solutions $\hat{x}$ from the robust-optimal solution set and their corresponding selection problem.}
	\label{tab01}
\end{table}

Selection problem $\text{SP}_1$ attempts finding the solution $\hat{x}^D$ that is feasible for the deterministic problem \eqref{blp} when the values of the deterministic variables are set equal to the values of the deterministic variables associated with $\mathcal{U}^*$; $\hat{x}^D$ belongs to the robust-optimal solution set and satisfies the feasibility constraints with $\delta_j = 0, \forall j \in J$. This solution, if it exists, has the best objective function and is feasible with respect to the deterministic problem since it considers feasible solutions in the deterministic feasible region $X$.

Selection problem $\text{SP}_2$ attempts finding the solution $\hat{x}^R$ for the deterministic problem \eqref{blp} when the values of the deterministic variables are set equal to the values of the deterministic variables associated with $\mathcal{U}^*$. The feasible region is relaxed to allow some infeasibilities with respect to the deterministic problem. This solution, if it exists, has the best objective function and may be infeasible with respect to the deterministic problem since it considers solutions in the robust feasible region (i.e., $\mathcal{X}$ or $\mathcal{X}''$).

Selection problems $\text{SP}_3$ and $\text{SP}_4$ attempt to select robust solutions $\hat{x}^{UB}$ and $\hat{x}^{LB}$ yielding the highest and lowest objective functions for any solution in $\mathcal{U}^*$, respectively. The following propositions allow identifying these solutions in a given $\mathcal{U}^*$ in linear time.

\begin{proposition}\label{upperbound}
	Given a robust-optimal solution set $\mathcal{U}^*$, an upper bound robust solution $\hat{x}^{UB} \in \mathcal{U}^*$ satisfies $f(\hat{x}^{UB}) \geq f(y), \forall y \in \mathcal{U}^*$ and its values are set as follows:	
	
	\begin{align}\label{pessimistic}
		\hat{x}_i^{UB} = \left\{\begin{matrix} x_i^* & \text{ if } i\in C \\ 1, & \text{ if } c_i \geq 0, i \in U \\ 0, & \text{ if } c_i < 0, i \in U. \end{matrix} \right.
	\end{align}
	
\end{proposition}

\begin{proof}
	Consider a solution $y \in \mathcal{U}^*$ such that $y \neq \hat{x}^{UB}$. By expression \eqref{pessimistic}, if $c_i \geq 0$ for $i \in U$ in $\hat{x}_U^{UB}$, $\hat{x}_i^{UB} = 1$ and $c_i\hat{x}_i^{UB} = c_i \geq c_iy_i$; similarly, if $c_i < 0$ for $i \in U$, $\hat{x}_i^{UB} = 0$ and $c_i\hat{x}_i^{UB} = 0 \geq c_iy_i$. Given that $\hat{x}_C^{UB} = y_C = x_C^*$, then $f(\hat{x}^{UB}) = \sum_{x_i \in \hat{x}^{UB}}c_ix_i = \sum_{x_i \in x_C^*}c_ix_i + \sum_{x_i \in \hat{x}_U^{UB}}c_ix_i \geq \sum_{x_i \in x_C^*}c_ix_i + \sum_{y_i \in y_U}c_iy_i= \sum_{y_i \in y}c_iy_i= f(y)$. Therefore, $f(\hat{x}^{UB}) \geq f(y), \forall y \in \mathcal{U}^*$. 
\end{proof}

\begin{proposition}\label{lowerbound}
	Given a robust-optimal solution set $\mathcal{U}^*$, a lower bound robust solution $\hat{x}^{LB} \in \mathcal{U}^*$ satisfies $f(\hat{x}^{LB}) \leq f(y), \forall y \in \mathcal{U}^*$ and its values are set as follows:
	\begin{align}\label{optimistic}
		\hat{x}_i^{LB} = \left\{\begin{matrix} x_i^* & \text{ if } i\in C \\ 0, & \text{ if } c_i \geq 0, i \in U \\ 1, & \text{ if } c_i < 0, i \in U. \end{matrix} \right.
	\end{align}
	
	Then $f(\hat{x}^{LB}) \leq f(y), \forall y \in \mathcal{U}^*$.
\end{proposition}

\begin{proof}
	Proof of Proposition \ref{lowerbound} is similar to the proof of Proposition \ref{upperbound}.
\end{proof}

The performance of these selected solutions will be illustrated later in the numerical study in the context of the Knapsack Problem.
	\section{Experimental Results} \label{expsection}

This section presents applications of RBIU-$\delta$ and RBIU-$\delta$-CC in the context of the knapsack problem and presents experimental results exemplifying the nature of solutions that are robust with respect to implementation uncertainties.

\subsection{Problem Formulation}

The deterministic KP is formulated as the following BLP:
\begin{alignat}{3}
	\max & \sum_{i =1}^{n} c_i x_i \label{kpobj}\\
	\text{s.t.} & \sum_{i =1}^{n} a_i x_i \leq b \\
	& x_i \in \{0,1\}, \forall i \in I. \label{kpbin}
\end{alignat}

Assuming that there exists at least one uncertain variable, the robust KP under implementation uncertainty (RKP-$\delta$) is the following:
\begin{alignat}{3}
	\min \gamma(x) &  = \max_{y\in \mathcal{U}(x_C)}\left\{\sum_{i =1}^{n} -c_i y_i\right\} \label{rkpobj} \\
	\text{s.t.} & \max_{y \in \mathcal{U}(x_C)}\left\{\sum_{i =1}^{n} a_i y_i\right\} \leq b + \delta \\
	& y_i \in \{0,1\}, \forall i \in I. \label{rkpcons3}
\end{alignat}

Using the linear reformulation of the RBIU-$\delta$ shown in Section \ref{MIPRob}, the RKP-$\delta$ can be rewritten as an equivalent deterministic KP with a smaller number of decision variables as follows:
\begin{equation}\label{reducedrkp}
	\begin{split}
		\max & \sum_{i \in C}c_ix_i \\
		\text{s.t.} & \sum_{i \in C} a_ix_i \leq b' \\
		& x_i \in \left\{0,1\right\}, \forall i \in C;
	\end{split}
\end{equation}

\noindent where $b' = b + \delta - \sum_{i \in U }a_i$. The objective function value level is computed as $\gamma(x) = \sum_{i \in C}c_ix_i$.

The combined feasibility and cardinality-constrained robust formulation, RKP-$\delta$-CC, follows:
\begin{alignat}{3}
	\min & \qquad\gamma \label{rkpccobj}\\
	\text{s.t.} & \sum_{i \in C}-c_i x_i + \max_{\left \{S_0 : S_0 \subseteq U, |S_0| \leq \Gamma \right \}}\left\{ \sum_{i \in S_0} \left(\frac{-c_i + |-c_i|}{2}\right) + \sum_{i \in U \setminus S_0 }c_i x_i \right\} \leq \gamma(x) \\
	& \sum_{i \in C}a_{i} x_i + \max_{\left \{S_j : S_j \subseteq U, |S_j| \leq \Gamma \right \}}\left\{\sum_{i \in S_j} \left (\frac{a_{i} + |a_{i}|}{2} \right ) + \sum_{i \in U \setminus S_j}a_{i} x_i \right\}  \leq b + \delta, \forall j \in J\\
	&  x_i \in \{0,1\}, \forall i \in I. \label{rkpccbin}
\end{alignat}

The RKP-$\delta$-CC can be solved using the equivalent linear reformulation in Section \ref{cardinality}.

\subsection{Objectives of The Experimental Study And Performance Metrics}

The objective of the experimental study is to illustrate the application of the methodology in the context of the knapsack problem and evaluate the performance and behavior of the proposed robust models under implementation uncertainty in this context. The objectives of the experiments are:

\begin{enumerate}
	\item Study the performance of robust solutions under different levels of implementation uncertainty. The robust solutions include:
	
	\begin{enumerate}
		\item Robust KP with feasibility parameter in \eqref{reducedrkp}, RKP-$\delta$.
		
		\item Robust KP with combined feasibility and cardinality constraint in \eqref{rkpccobj}-\eqref{rkpccbin}, RKP-$\delta$-CC.
	\end{enumerate}

	\item Study the performance of the different selected prescribed solutions in Table \ref{tab01}.
	
	\item Study of the probability upper-bound of the RKP-$\delta$-CC derived in Section \ref{bounds}.
\end{enumerate}

The performance of the deterministic solution will be reported as a baseline.

Performance is measured in terms of the objective function and feasibility. To measure the aggregate objective function performance of a robust solution set, we define the average objective performance ratio as follows; for a given solution set $\mathcal{U}(x)$:
\begin{equation}\label{optlevelexp}
	\bar{\Delta}(\mathcal{U}(x))= 1 - \frac{\sum_{y \in \mathcal{U}(x)}f(y)}{f(x^*_{KP})|\mathcal{U}(x)|};
\end{equation}

\noindent where $x^*_{KP}$ is the solution of the deterministic KP in \eqref{kpobj}-\eqref{kpbin}. The smaller $\bar{\Delta}$ is, the better the performance of the robust solution set.

Given a solution set $\mathcal{U}(x)$, the feasibility level $h(\mathcal{U}(x))$ is defined as the proportion of the solutions in $\mathcal{U}(x)$ that is feasible with respect to the deterministic feasible set as follows:
\begin{equation}\label{feaslevelexp}
	h(\mathcal{U}(x)) = \frac{|\mathcal{U}(x) \cap X|}{|\mathcal{U}(x)|}.
\end{equation}

To evaluate the performance of a single prescribed solution, $x$, instead of the aggregate performance of a robust set, the objective performance ratio with respect to the deterministic solution:
\begin{equation}
	\Delta(x) = 1 - \frac{f(x)}{f(x^*_{KP})}.
\end{equation}

We also evaluate the performance of a single solution, $x$, by measuring the level of the infeasibility of a robust solution set in terms of the level of constraint violation as follows:
\begin{equation}
	F(x) = \max\left\{\frac{\sum_{i = 1}^{n}a_{i}x_{i} - b}{b},0 \right\}.
\end{equation}

\subsection{Simulation of Uncertainty And Test Problem Generation}

The level of uncertainty during the modeling of the robust formulations is defined by the number of uncertain variables $|U|$; the feasibility level of the RKP-$\delta$ by $\delta$; the level of impact of uncertainty and feasibility of the RKP-$\delta$-CC by $\Gamma$ and $\delta$, respectively. The experimental statistic $p=\Gamma/|U|$ measures the level of conservatism of the RKP-$\delta$-CC's solutions by defining the percentage of uncertain variables that may be affected by implementation uncertainty.

$|U_S|$ is the set of simulated uncertain variables. In this experiment, full enumeration is considered to evaluate the solutions' performances for values of $|U_S| \leq 10$; otherwise, a random sampling of size $2^{10}$ is used to evaluate the solutions' performances to reduce the computational time of enumerating each element in $\mathcal{U}(x)$. While sampling, the implemented value of an uncertain variable may be equal or different to the prescribed value with the same probability of 0.5. The description of the simulation and corresponding pseudo-code describing the simulation is shown in the Appendix.

This experiment consists of twelve replications; each replication consists of a KP instance with $n = 100$. The costs $c_i$ are randomly chosen from the set $\left\{21,22,...,80\right\}$, the weights $a_i$ are randomly chosen from the set $\left\{41,42,...,60\right\}$, and the capacity $b$ is equal to $0.5 \sum_{i=1}^{n}a_i$.

The RKP-$\delta$ and RKP-$\delta$-CC are formulated considering $u = |U|/n = 1\%,3\%, 5\%, 7\%, 9\%$ of the total number of variables $n$. The uncertain variables are selected from 1 to $|U|$, and deterministic variables from $(|U|+1)$ to $n$. The feasibility parameter $\delta$ is set at $0\%,1\%,2\%,3\%,4\%,5\%$ of the value of $b$. The values of $p$ considered are 40\%, 60\% and 80\%.

\subsection{Experimental Analysis}


\subsubsection{RKP-$\delta$ With Feasibility Performance}\label{rkpfeas}

\begin{figure}[h]
	\centering \includegraphics[scale = 0.64]{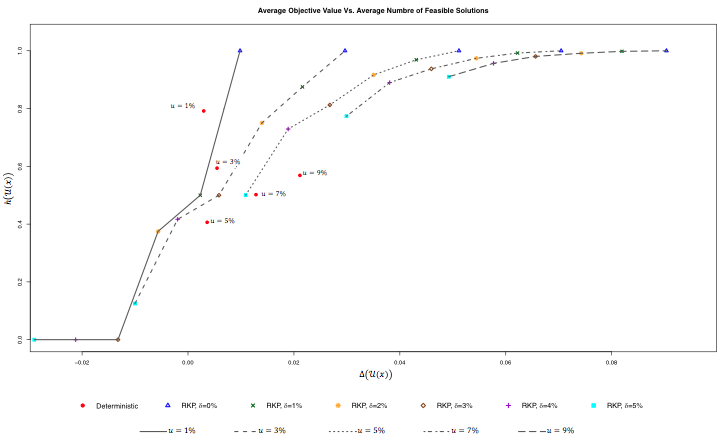}
	\caption{RBK-$\delta$ for uncertainty level $u$ and feasibility $\delta$ combinations. \label{plot03}}{}
\end{figure}

We run an experiment to evaluate the performance of RPK-$\delta$ \eqref{rkpobj}-\eqref{rkpcons3} using different levels of different levels of uncertainty $u$ and feasibility parameter $\delta$. Figure \ref{plot03} summarizes the results in a plot of feasibility level \eqref{feaslevelexp} versus average objective performance \eqref{optlevelexp}; in this plot, markers represent different $\delta$s, and lines represent different uncertainty levels. The average performance of the deterministic solutions are shown as solid circles for different uncertainty levels.

All results using $\delta=0$ have perfect feasibility performance for all levels of uncertainty tested, confirming that RBK-$\delta$ can guarantee deterministic-feasible outcomes. The feasibility performance decreases as $\delta$ increases; this effect accentuates as the uncertainty level decreases. For the data tested, most runs with $u>3\%$ yielded $80\%$ or better feasibility levels than the deterministic solutions that could achieve at most $60\%$ feasibility; furthermore, the objective function for these runs ranged between 2\%-10\% from the deterministic solution.

The losses in feasibility caused by increases in $\delta$ diminish as the uncertainty levels increase. This suggests that as more variables are uncertain, there are more opportunities to find robust solutions, especially as $\delta$ expands the robust feasible solution $\mathcal{X}$ beyond the deterministic feasible region $X$. For instance, results for $u=9\%$ improve their objective performance from about $9\%$ down to about $5\%$ of the optimal objective value, when $\delta$ increases from 0 up to 5\%, while only decreasing the feasibility level from 1.0 down to 0.9. This demonstrates the practical potential of RBK-$\delta$ to generate robust solutions.

The negative effects of $\delta$ are magnified for low uncertainty levels; i.e., for $u<3\%$, RBK-$\delta$ can rapidly produce over-optimistic infeasible solutions with increasing $\delta$s. Moreover, over-optimistic solutions display worst feasibility performance than deterministic solutions. This is expected in knapsack problems if the number of variables that can flip values is small because this leads to fewer outcome combinations reducing the opportunities to regain feasibility. Hence, $\delta$s must be increased carefully (e.g., in finer increments) at low levels of uncertainty.

In summary, RBK-$\delta$ can produce practical robust solutions with acceptable objective performance and superior feasibility performance compared to the deterministic solution. Results also suggest that the effect of the max-infeasibility parameter $\delta$ enables finding better robust solutions by considering robust regions beyond the deterministic feasible region.

\subsubsection{RKP-$\delta$-CC And Interactions With RKP-$\delta$}\label{exps2}

\begin{figure}[h]
	\centering \includegraphics[scale = 0.58]{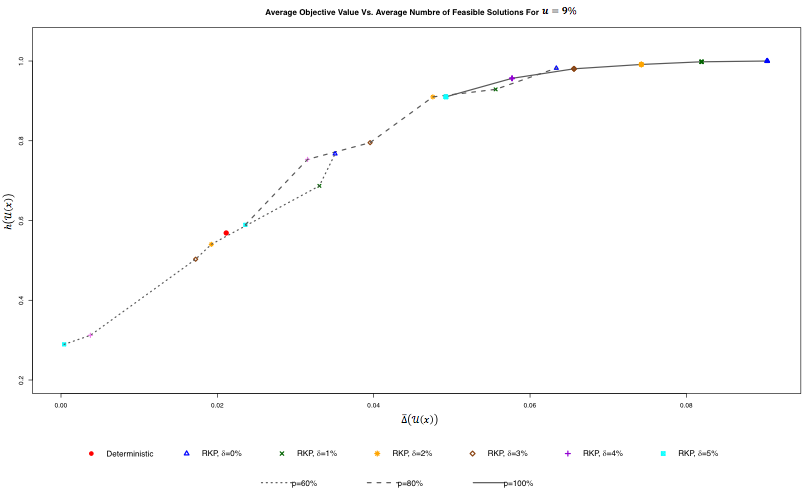}
	\caption{RKP-$\delta$-CC, $u=9\%$, and combinations of degree of conservatism $p$ and feasibility $\delta$. \label{plot04}}{}
\end{figure}

Figure \ref{plot04} depicts the performance of RKP-$\delta$-CC in the feasibility level versus objective performance plane for different levels of uncertainty and degree of conservatism, $p = \Gamma/|U|$, for the same feasibility level $\delta=0$\%. Notice that when $p = 100\%$, RKP-CC performs the same as RKP-$\delta$, the performance of the deterministic and RKP-$\delta$ solutions are also shown for reference.

We run an experiment to evaluate the performance of RPK-$\delta$-CC \eqref{rkpccobj}-\eqref{rkpccbin} for $u=9\%$ different combinations of feasibility parameter $\delta$ and level of conservatism $p=\Gamma/|U|$. Figure \ref{plot04} summarizes the results in a plot of feasibility level \eqref{feaslevelexp} versus average objective performance \eqref{optlevelexp}; in this plot, markers represent different $\delta$s, and lines represent different levels of conservatism $p$. The average performance of the deterministic solutions are shown as a solid circle.

In this figure the curve for $p=1.0$ correspond to the curve for $\delta=0$ in Figure \ref{plot03}. In Figure \ref{plot04} this curve illustrates how decreasing the level of conservatism (i.e., the cardinality-constrained parameter) yields better objective function performance while introducing infeasibilities for runs with minimum $\delta=0$. Likewise, one could also state that given maximum conservatism (i.e., $p=1.0$), incrementing the feasibility $\delta$ also causes better objective function performance while introducing infeasibilities. In fact, for a given uncertainty level, observing how the curves for different $p$ overlap with each other suggest that different combinations $(p,\delta)$ produce similar solutions. This observation suggests that both the max-infeasibility and the cardinality-constrained approaches have a similar effect in generating less conservative solutions. We have observed the same phenomenon occurring at $u=3\%$ and $u= 5\%$.

By combining both parameters $p$ and $\delta$, we provide the decision-maker with the possibility of defining a degree of conservatism in the sense of \cite{bertsimas2004price} by specifying a given $p$, while using $\delta$ to expand the robust feasible space $\mathcal{X}$ beyond the deterministic feasible region $X$ to find robust solutions with even better objective performance with only small degradation of the feasibility performance. Consider, for instance, a decision-maker with a degree of conservatism $p=0.8$; without the effect of $\delta$, this CC solution will have an objective performance of about 0.065 from the optimal solution and excellent feasibility performance 1.0. Increasing $\delta$ from 0 to 2\% would improve the objective performance from 0.065 to less than 0.05 while only reducing the feasibility performance down to above 0.9. Notice that this is significantly better than the deterministic solution performance with less than 0.6 feasibility level. 

In summary, using RKP-$\delta$-CC provides additional control to find robust solutions that are more practical and improve those that would find the cardinality-constrained approaches.

\subsubsection{Performance of Selected Prescribed Solutions}\label{exps3}

The experiments in Sections \ref{rkpfeas} and \ref{exps2} evaluate the aggregate expected performance of the robust solution set $\mathcal{U}^*$. In this section, we demonstrate the advantages of selecting solutions as an optimization over the robust set by evaluating the performance of selected member solutions described in Section \ref{solselection}. Specifically, selected solutions from $\mathcal{U}^*_{\delta}$ and $\mathcal{U}^*_{\delta-CC}$ associated with formulations considering feasibility \eqref{reducedrkp}, and combined feasibility and cardinality-constrained \eqref{rkpccobj}-\eqref{rkpccbin}, respectively.

\begin{table}[h]
	\centering
	\begin{tabular}{| c | p{8cm} |}
		\hline
		Prescribed solution $\hat{x}$ & Selection Problem, SP \\ \hline
		$\hat{x}^D_{\delta}$ & $\text{SP}_1$: Solve \eqref{blp}; set $x_C \leftarrow x_C^*; (\hat{x}^*_C,\hat{x}^*_U)\in \mathcal{U}^*_{\delta}$ \\
		$\hat{x}^R_{\delta}$ & $\text{SP}_2$: Solve \eqref{blp}; set $x_C \leftarrow x_C^*$  and the right-hand side $\leftarrow b + \delta; (\hat{x}^*_C,\hat{x}^*_U)\in \mathcal{U}^*_{\delta}$\\
		$\hat{x}^{UB}_{\delta}$ & $\text{SP}_3$: Find $x\in \mathcal{U}^*_{\delta}:f(x) \geq f(y), \forall y \in \mathcal{U}^*_{\delta}$ \\
		$\hat{x}^{LB}_{\delta}$ & $\text{SP}_4$: Find $x\in \mathcal{U}^*_{\delta}:f(x) \leq f(y), \forall y \in \mathcal{U}^*_{\delta}$ \\
		$\hat{x}^D_{\delta-CC}$ & $\text{SP}_5$: Solve \eqref{blp}; set $x_C \leftarrow x_C^*; (\hat{x}^*_C,\hat{x}^*_U)\in \mathcal{U}^*_{\delta-CC}$ \\
		$\hat{x}^R_{\delta-CC}$ & $\text{SP}_6$: Solve \eqref{blp}; set $x_C \leftarrow x_C^*$  and the right-hand side $\leftarrow b + \delta; (\hat{x}^*_C,\hat{x}^*_U)\in \mathcal{U}^*_{\delta-CC}$\\
		$\hat{x}^{UB}_{\delta-CC}$ & $\text{SP}_7$: Find $x\in \mathcal{U}^*_{\delta-CC}:f(x) \geq f(y), \forall y \in \mathcal{U}^*_{\delta-CC}$ \\
		$\hat{x}^{LB}_{\delta-CC}$ & $\text{SP}_8$: Find $x\in \mathcal{U}^*_{\delta-CC}:f(x) \leq f(y), \forall y \in \mathcal{U}^*_{\delta-CC}$ \\ \hline		
	\end{tabular}
	\caption{Selected prescribed solutions $\hat{x}$ from the robust-optimal solution set and their corresponding selection problem.}
	\label{tab02}
\end{table}

Given a selected prescribed solution $\hat{x}$, we compute the corresponding feasibility and objective performances as prescribed; in turn, the expected performance in the face of implementation uncertainty $E[\tilde{x}]$ is estimated using simulation to observe the behavior of its corresponding implemented solution $\tilde{x}$. Finally, the behavior of the deterministic KP solution is also evaluated as a baseline.

\begin{figure}[h]
	\centering \includegraphics[scale = 0.49]{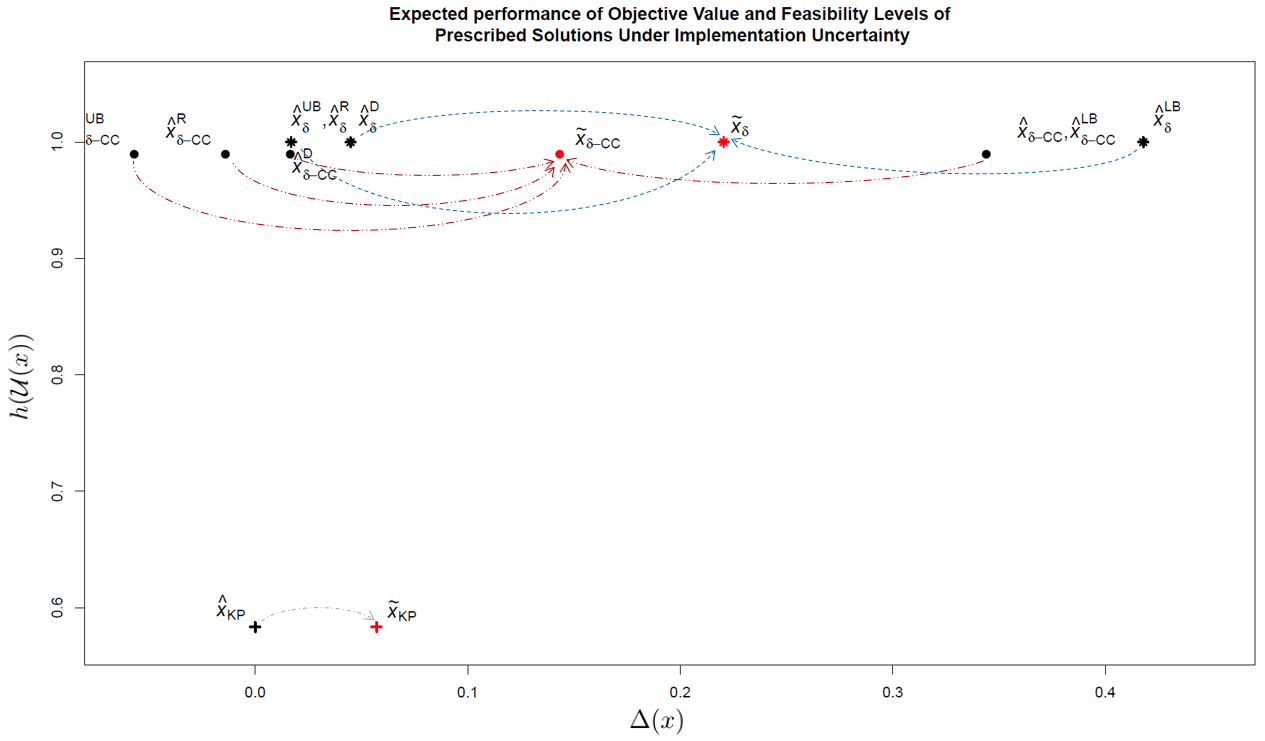}
	\caption{$h(\mathcal{U}(x))$ versus $\Delta(\mathcal{U}(x))$. Prescribed versus implemented performance of selected robust solutions. Plot for $u=25\%, \delta = 10\%$ and $p = 80\%$. \label{plot07}}{}
\end{figure}

Figure \ref{plot07} summarizes the results in the $h(\mathcal{U}(x))$ versus $\Delta(\mathcal{U}(x))$ plane for $u=25\%, \delta = 10\%$ and $p = 80\%$. As expected, all selected solutions from the same robust solution set converge to the aggregate solution of the robust solution set; also, as described in the previous experiments, all robust solutions for these levels of uncertainty, feasibility, and cardinality-constrained parameters have either perfect or close-to-perfect feasibility with respect to the deterministic KP. With relatively good objective value performance, the prescribed and implemented deterministic solution is feasible only in 0.60 cases when affected by implementation uncertainty.

\begin{figure}[h]
	\centering \includegraphics[scale = 0.49]{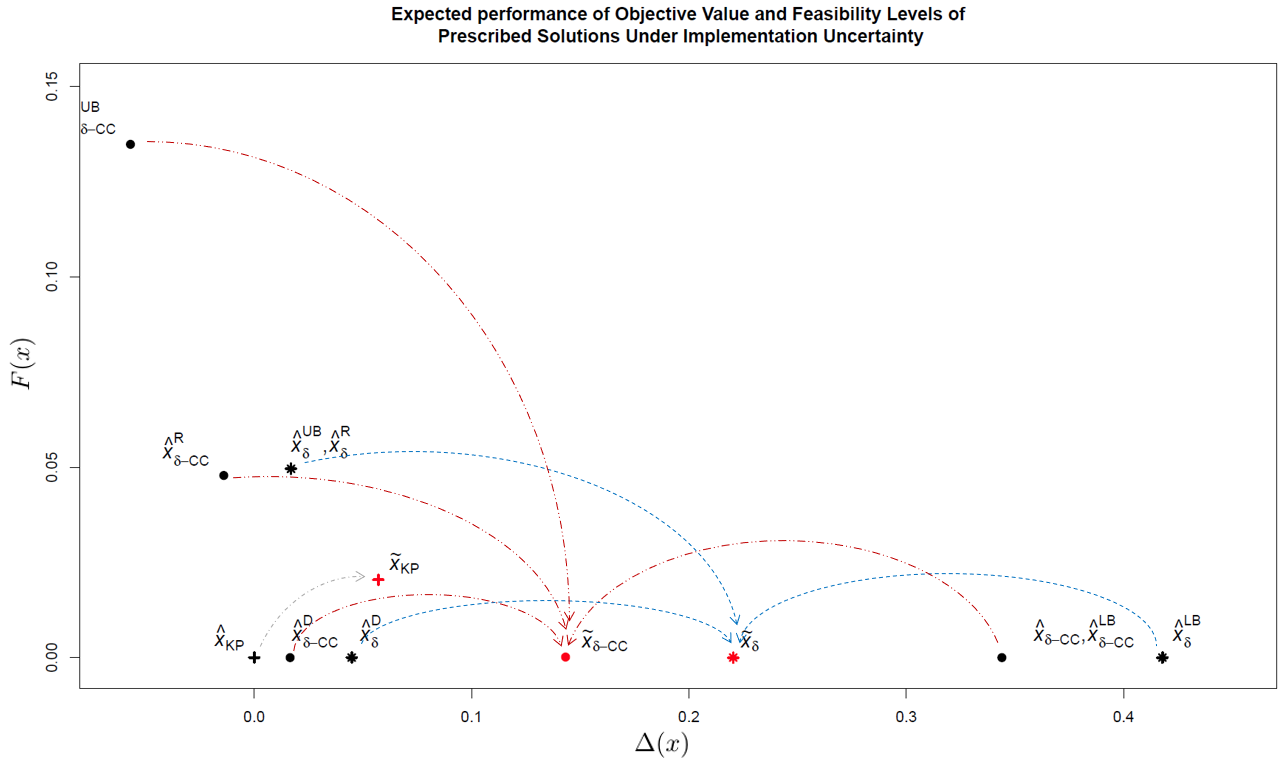}
	\caption{$F(x)$ versus $\Delta(\mathcal{U}(x))$. Prescribed versus implemented performance of selected robust solutions. Plot for $u=25\%, \delta = 10\%$ and $p = 80\%$. \label{plot08}}{}
\end{figure}

Figure \ref{plot08} summarizes the result in the $F(x)$ versus $\Delta(\mathcal{U}(x))$ plane for the same conditions of $u$, $\delta$ and $p$. In this plot, one can observe the magnitude of feasibility violation of the robust solutions as prescribed. For instance, although $\hat{x}^R_{\delta-CC}$ has a prescribed value that seems better to the ideal deterministic value, it does so by compromising feasibility with $F(\hat{x}^R_{\delta-CC}) > 0$; this may discourage a decision-maker from implementing it even though its expected performance (Figure \ref{plot07}) will be almost always feasible with respect to the deterministic problem. Alternatively, in this case, $\hat{x}^D_{\delta-CC}$ has a similar prescribed objective value and feasibility performance as the ideal deterministic solution. However, in the face of implementation uncertainty, $\hat{x}^D_{\delta-CC}$ will guarantee feasibility and maintain its objective performance in the best case, and only deteriorate its objective value down to $E\left[\hat{x}^D_{\delta-CC}\right]$ on the average; hence, $\hat{x}^D_{\delta-CC}$ might be a good candidate for implementation when feasibility and objective performance are both important implementations considerations. Furthermore, in this case, if some degree of violation was acceptable in a particular implementation, solutions $\hat{x}^{UB}_{\delta-CC}$ and $\hat{x}^R_{\delta-CC}$ may be of interest for implementation since in the best case, they may result in better objective performance than the ideal deterministic problem.

In summary, experimental results support the advantages of selecting solutions by optimization over the robust set. They suggest that selected solutions possess the desired robustness qualities and desirable implementation characteristics such as similar or better performance as the ideal deterministic solution.

\subsection{RKP-$\delta$-CC Upper Bound Probability}

Table \ref{tab09} displays a plot of the results of the test conducted to validate the quality of the probability upper bounds developed in Section \ref{bounds}.

Experimental results suggest that the theoretical bounds \eqref{ubprob} of the probability that RKP-$\delta$-CC solutions become infeasible are accurate for high values of $p$; however, they may be less accurate when $p$ is small. Table \ref{tab09} presents the differences between the theoretical and experimental probabilities,$\Delta P(\text{infeasible})$, that RKP-$\delta$-CC solutions become infeasible for different values of $p$ and $\delta = 0$. Results in Table \ref{tab09} show that the difference is large for $p = 40\%$ and decreases as the value of $p$ increases; for $p = 80\%$, there does not exist a significant difference between the theoretical and experimental probabilities.

\begin{table}[h]
	\centering
	\begin{tabular}{| c | c | c | c |}
		\hline
		$u$ & $p = 0.4$ & $p = 0.6$ & $p = 0.8$ \\ \hline
		5\% & 0.125 & 0.015625 & 0.002604167 \\
		15\% & 0.1824646 & 0.019856771 & 0.000518799 \\
		25\% & 0.103495717 & 0.03997306 & 0.000129739 \\
		35\% & 0.119405874 & 0.043705932 & 5.84209E-05 \\
		45\% & 0.122585696 & 0.041373798 & 7.68705E-06 \\ \hline
	\end{tabular}
	\caption{Differences between the theoretical and experimental probabilities that RKP-$\delta$-CC solutions become infeasible.}
	\label{tab09}
\end{table}

	\section{Conclusions}\label{conclu}

The prescribed solutions of optimization problems may require changes during their implementation due to unforeseen reasons not considered during the optimization model development. We model possible changes of a prescribed solution as implementation uncertainty affecting the decision variables; this paper focuses on RBI-$\delta$U that limits the scope of the problem to the context of binary linear programs, and the case where there is no information about the probability distributions describing the uncertainties. The main solution approach is robust optimization with a distinct characteristic that uncertainty affects the decision variables rather than the model parameters. The solution methodology allows controlling the level of conservativism associated with worst-case solutions; specifically, experimental results in the context of the knapsack problem suggest that robust solutions display beneficial properties when compared to the deterministic counterparts in terms of both feasibility and objective function robustness when facing implementation uncertainties. The experiments also suggest the methodology is not appropriate under all problem conditions, but that is clearly superior to the deterministic solutions as the importance of feasibility of the solution at implementation and the degree of uncertainty increase. The experiments also suggest that the mechanisms to control conservativism, namely the feasibility and cardinality-constrained parameters provided are effective to broaden the applicability region of the methodology.

The robust solution set may be large as it grows exponentially with the number of uncertain variables. Therefore, we propose selecting specific solutions in the robust solution set such that the prescribed (selected) solution also possesses desirable characteristics for implementation. Experimental results suggest that the robust set may contain robust solutions that have a close-to-optimal performance with respect to the deterministic version of the problem.

This work opens opportunities for future research; for instance, the development of measures of objective degradation less conservative than the worst-case objective value. Of interest is also solving the optimization problem over the robust set directly, i.e. without having to identify the robust set as an intermediate step. Also, there exists an opportunity for the application of the proposed methodologies in other application contexts.


\section*{Appendix}\label{simulation}

	When simulating uncertainty, a full enumeration is used for $|U_S| \leq 10$ and random sampling otherwise. Algorithm \ref{simulalgo} presents the pseudo-code describing the simulation of uncertainty.

\begin{breakablealgorithm}
		\caption{Simulation of uncertainty}
		\label{simulalgo}
		\textbf{input:} a KP instance, a solution $x=(x_C,x_U)$ and the simulated level of uncertainty $|U_S|$. \\
		\textbf{output:} Aggregated objective performance $\bar{\Delta}(\mathcal{U}(x))$ and feasibility level $h(\mathcal{U}(x))$.
		\begin{algorithmic}[1]
			\Function {Simulation of uncertainty}{Instance of KP, $x$, $|U_S|$}
			\If {$|U_S|\leq 10$} \Comment{Full enumeration is performed if $|U_S|\leq 10$}
			\ForEach {$y \in \mathcal{U}(x)$}
			\State Compute $\sum_{i = 1}^{n}a_iy_i$
			\If {$\sum_{i = 1}^{n}a_iy_i \leq b$}
			\State Add 1 to \textit{Feasible\_Solutions}
			\EndIf
			\State Compute $f(y) = \sum_{i = 1}^{n}c_iy_i$
			\State Add $f(y)$ to \textit{Total\_Objective}
			\EndFor
			\State $h(\mathcal{U}(x))=\frac{\textit{Feasible\_Solutions}}{2^{|U_S|}}$
			\State $\bar{\Delta}(\mathcal{U}(x))=f(x_D^*)-\frac{\textit{Total\_Objective}}{2^{|U_S|}}$
			\Else \Comment{Random sampling is performed if $|U_S|> 10$}
			\For {1 to $2^{10}$} \Comment{The size of the sample is $2^{10}$}
			\State Define $y = (x_C,y_U)$ \Comment{Sample vector with the deterministic vector equal to the prescribed solution}
			\ForEach {$y_i$ in $y_u$}
			\If {random value between in $[0,1] > 0.5$} \Comment{Uncertain variables change value with probability 0.5 } 
				\State $y_i = x_i, x_i$ in $x_U$ \Comment{Uncertain variable does not change its value}
			\Else
				\State $y_i = 1-x_i, x_i$ in $x_U$ \Comment{Uncertain variable changes its value}
			\EndIf
			\EndFor						
			\State Compute $\sum_{i = 1}^{n}a_iy_i$
			\If {$\sum_{i = 1}^{n}a_iy_i \leq b$}
			\State Add 1 to \textit{Feasible\_Solutions}
			\EndIf
			\State Compute $f(y) = \sum_{i = 1}^{n}c_iy_i$
			\State Add $f(y)$ to \textit{Total\_Objective}
			\EndFor
			\State $h(\mathcal{U}(x))=\frac{\textit{Feasible\_Solutions}}{2^{10}}$
			\State $\bar{\Delta}(\mathcal{U}(x))=f(x_D^*)-\frac{\textit{Total\_Objective}}{2^{10}}$
			\EndIf
			\State \textbf{return} $\bar{\Delta}(\mathcal{U}(x))$ and $h(\mathcal{U}(x))$
			\EndFunction
		\end{algorithmic}
\end{breakablealgorithm}

	\bibliographystyle{informs2014} 
	\bibliography{references} 

\end{document}